\documentclass[a4paper,10pt,leqno,english]{smfart}
\usepackage{amsmath,a4wide}
\usepackage{amssymb}
\usepackage{mathrsfs}
\usepackage{inputenc,xspace}
\usepackage{euscript}
\usepackage[all]{xy}
\usepackage[english,frenchb]{babel}
\usepackage{graphicx}
\usepackage[T1]{fontenc}
\usepackage{amscd}
\usepackage[all]{xy}
\usepackage{graphicx}

\newtheorem{prop}{Proposition}
\newtheorem{theo}[prop]{Theorem}
\newtheorem{lemm}[prop]{Lemma}
\newtheorem{coro}[prop]{Corollary}
\theoremstyle{definition}

\newtheorem{rem}[prop]{Remark}
\newtheorem{exa}[prop]{Example}

\def\CC{\bold C}
\def\P{\bold P^1}
\def\QQ{\overline{\bold Q}}
\def\PP{\bold P^*}
\def\Ga{\Bbb G_{\rm a}}
\def\tors{{\rm tors}}
\def\Div{{\rm Div}}
\def\GQ{{G_{\bold Q}}}

\begin{document}
\title[Arithmetic properties of Lam\'e operators]{{\small\it\bfseries Secondo convegno italiano di teoria dei numeri\\ Parma, 13-15th of  november 2003\\ \vskip.5cm\hrule\vskip.7cm}
{\larger Some arithmetic properties of Lam\'e operators with dihedral monodromy}}
\author{Leonardo Zapponi}
\address{EPFL, Lausanne}
\email{leonardo.zapponi@epfl.ch}
\urladdr{http://alg-geo.epfl.ch/~leo/leonardo.html}
\date{\today}
\begin{abstract} In this paper, we describe some arithmetic properties of Lam\'e operators with finite dihedral projective monodromy. We take advantage of the deep link with Grothendieck's theory of dessins d'enfants, following~\cite{Litcanu_1,Litcanu_2}. We focus more particularly on the case of projective monodromy of order $2p$, where $p$ is an odd prime number.
\end{abstract}
\maketitle

\section*{Introduction}
Lam\'e operators are a particular class of second order Fucshian differential operators on the projective line. In some special cases, they admit a complete system of algebraic solutions, i.e. they have a finite monodromy. This last question has been intensively studied by F. Baldassarri, B. Chiarellotto, B. Dwork and more recently by F. Beukers, S. Dahmen, R. Li\c tcanu, A. van der Waall. It turns out that there are finitely many equivalence classes of operators with fixed finite monodromy group and that these are automatically defined over a number field. In particular, there is a well defined action of the absolute Galois group on these objects. The enumeration of Lam\'e operators with finite (projective) monodromy has been one of the main motivations in this topic. Recently, a deep link with Grothendieck's theory of dessins d'enfants appeared; this point of view has been succesfully adopted by R. Li\c tcanu in~\cite{Litcanu_1,Litcanu_2}, allowing an explicit and combinatorial enumeration (see also~\cite{Dahmen}).

In this paper, we focus on the case of Lam\'e operators with exponent $n=1$ having finite dihedral projective monodromy group. In \S1 we give some basic definitions, by introducing the field of moduli of a Lam\'e operator, which is the smallest field of definition and is invariant under equivalence. In \S2, we translate a criterion of F. Baldassarri~\cite{Baldassarri} in terms of generalized jacobians. This criterion asserts that the existence of a Lam\'e operator with dihedral monodromy of order $2N$ is related to the existence of a $2N$-torsion point on an elliptic curve with some extra properties. In \S3, inspired by the work of R. Li\c tcanu, we briefly descibe the correspondence between the set of Lam\'e operators with dihedral monodromy and a particular class of dessins d'enfants. We then give two direct applications, namely the finiteness of the set of equivalence classes of Lam\'e operators with fixed projective monodromy and the fact that the fields of moduli of such operators are number fields. In \S4 we prove, following a celebrated result of L. Merel~\cite{Merel}, that there exist finitely many equivalence classes of Lam\'e operators with dihedral monodromy and field of moduli of bounded degree. We then investigate more closely the behaviour of the field of moduli, by taking advantage of some recent developments on the study of (semi-stable models of) covers between curves. First of all, a result of S. Beckmann~\cite{Beckmann} (see also~\cite{Chambert-Loir, Emsalem, Wewers}) implies that it is unramified outside the primes which are less than or equal to half the order of the (dihedral) monodromy group. In \S5 we study the case of dihedral monodromy of order $2p$, where $p$ is an odd prime number. Fist of all, we show that in this case the field of moduli is effectively ramified at the primes lying above $p$ (by giving a lower bound for their ramification index). We then prove that the elliptic curve which is naturally attached to the operator always has potentially good reduction at these primes and we give a supersigularity criterion for the reduced curve. The first of these last results follow from~\cite{Zapponi_1} but it can also easily be deduced from~\cite{Wewers}, while the last criterion need a more accurate investigation on the action of the Cartier operator (which is done in~\cite{Zapponi_3}). Finally, we show that the elliptic curve admits a smooth model at a prime $\frak p$ of the field of moduli lying above $p$ if and only if the ramification index of $\frak p$ is large enough. We then illustrate the results with the description of the Lam\'e operators with dihedral projective monodromy of order $14$.

This paper is the result of a work which is still in progress: for example, it is now possible to completely determine the ramification index of the primes in the field of moduli lying above $p$ (in the case of dihedral monodromy of order $2p$); these results being not yet published, we decided not to include them.

I would like to thank the organizers and the participants of the {\it Secondo convegno italiano di teoria di Numeri}, yeld in Parma during the month of November 2003. A special thank goes to A. Zaccagnini for his welcome and his local organization and to R. Li\c tcanu for the instructive discussions and comments on the subject.

\section{Lam\'e operators and their fields of moduli}
 A {\it Lam\'e operator} is a second order differential operator on the projective line defined by
$$L_n=L_{n,g_2,g_3,B}=D^2+\frac{f'}{2f}D-\frac{n(n+1)x+B}f$$
where $D=d/dx$, $f(x)=4x^3-g_2x-g_3\in\CC[x]$ with $\Delta=g_2^3-27g_3^2\neq0$ and $B\in\CC$. Let $E$ be the elliptic curve defined by the affine equation $y^2=f(x)$, denote by $0_E$ its origin (the point at infinity) and by $\sigma$ the canonical involution $\sigma(x,y)=(x,-y)$. We say that the operator $L_n$ is {\it associated} to $E$ and that $B$ is the {\it accessory parameter}. Two Lam\'e operators $L_{n,g_2,g_3,B}$ and $L_{n,g_2',g_3',B'}$ are {\it equivalent} (or {\it scalar equivalent}, following~\cite{Beukers}) if there exists $u\in\CC^*$ such that $g_2'=u^2g_2$, $g_3'=u^3g_3$ and $B'=uB$. In This paper, we are mainly concerned with the case $n=1$ but many results can be carried over to the general case.

\begin{lemm}\label{lemm1} There is a natural bijection between the equivalence classes of Lam\'e operators $L_1$ and the isomorphism classes of pairs $(E,P)$ where $E$ is an elliptic curve and $P\neq0_E$ is a point on it.
\end{lemm}

\begin{proof} Given a pair $(E,P)$, we may suppose (since we are working up to isomorphism) that $E$ is given by a Weierstrass model $y^2=4x^3-g_2x-g_3$; we then associate to it the Lam\'e operator $L_{1,g_2,g_3,x(P)}$. Conversely, given $L_1=L_{1,g_2,g_3,B}$, we consider the pair $(E,P)$ where $E$ is the elliptic curve associated to $L_1$ and $P$ is one of the two points of $E$ for which $x(P)=B$. One easily checks that equivalent Lam\'e operators correspond to isomorphic pairs and vice-versa.
\end{proof}

The field $\bold Q(g_2,g_3,B)$ is the {\it field of definition} of the operator $L_1$, its {\it field of moduli} $K$ is the intersection of the fields of definition of all the Lam\'e operators equivalent to it; it contains the field $\bold Q(j)$, where $j=1728g_2^3/\Delta$ is the absolute modular invariant associated to the elliptic curve and one can easily prove that $K$ is a actually a field of definition. More explicitely, we find $K=\bold Q(j,j_1,j_2,j_3)$, where we have set $j_1=B^4g_2/\Delta, j_2=B^2g_2^2/\Delta$ and $j_3=B^3g_3/\Delta$. It is possible to define the field of moduli of a pair $(E,P)$ which coincides with the field of moduli of the corresponding Lam\'e operator (following Lemma~\ref{lemm1}).

\begin{rem} The above Lemma asserts that equivalence classes of Lam\'e operators bijectively corresponds to the $\CC$-rational points of the moduli space $\mathcal M_{1,2}$ (the marked point $0_E$ is implicitely given in the definition of $E$). The field of moduli of an operator is just the field of definition of the corresponding point. 

\end{rem}

\section{dihedral projective monodromy and generalized Jacobians}

Let $E$ be an elliptic curve defined by a Weierstra{\ss} equation, as in~\S1. Recall that the {\it generalized Jacobian} $J_\frak m$ associated to the modulus $\frak m=2[0_E]$ is the quotient of the group of degree zero divisors of $E$ which are prime to $0_E$ with respect to the group of principal divisors of the type $D=(t)$ with $t$ regular at $0_E$ and $v_{0_E}(dt)\geq1$ (we refer to~\cite{Serre} for a detailed exposition on this subject). In particular, we have an exact sequence of algebraic groups
$$0\to\Ga\to J_\frak m\stackrel\pi\longrightarrow E\to0$$
There is a natural map $\varphi:E\to J_\frak m$ which sends a point $P$ to the equivalence class of the divisor $[P]-[\sigma(P)]$. It is important to note that even if the composition $\pi\circ\varphi$ is the multiplication by $2$ map, the morphism $\varphi$ is not a homomorphism between algebraic groups. The following result is a reformulation of the existence criterion in~\cite{Baldassarri} for operators $L_1$ with dihedral projective monodromy group. As usual, given a group $G$, we denote by $G[n]$ its $n$-torsion subgroup.

\begin{prop}\label{prop1} Let $E$ be an elliptic curve and $P$ a point on it. The following conditions are equivalent:
\begin{itemize}
\item The Lam\'e operator associated to the pair $(E,P)$ (cf. Lemma~\ref{lemm1}) has dihedral projective monodromy of order $2N$.
\item The element $\varphi(P)$ has exact order $N$.
\end{itemize}
In particular, if one of these condition is fulfilled then $P$ is a $2N$-torsion point on $E$.\end{prop}

\begin{proof} We know from~\cite{Baldassarri} that the operator $L_1$ attached to $(E,P)$ has dihedral projective monodromy of order $2N$ if and only if $P\in E[2N]\setminus E[2]$ satisfies the following conditions:
\begin{enumerate}
\item The point $Q=2P\in E[N]$ has exact order $N$.
\item Setting $D=[P]-[\sigma(P)]$, we have $ND=(t)$ with $v_{0_E}(dt)\geq2$.
\end{enumerate}
In terms of generalized Jacobians, these two conditions can be restated by saying that $D$ defines a point of exact order $N$ in $J_\frak n$, with $\frak n=3[0_E]$. For any $P\in E[2N]$, there exists a unique function $t$ for which $ND=(t)$ and $t(0_E)=1$. Since $\sigma^*D=-D$, we obtain $\sigma^*t=1/t$. In particular, setting $\omega=dt/t$, we have $\sigma^*\omega=-\omega$, so that, if $z=x/y$ is the usual uniformizer at $0_E$, we get the formal expansion
$$\omega=(a_0+a_2z^2+a_4z^4+\dots)dz$$
The condition $v_{0_E}(dt)\geq2$ can be restated as $v_{0_E}(dt/t)\geq2$ and the above expression implies that it is equivalent to $v_{0_E}(dt/t)\geq1$, as desired.
\end{proof}

 Fix an element $\tau$ of the upper half plane corresponding to $E$, so that the elliptic curve is isomorphic to the quotient $\bold C/\Lambda$, where $\Lambda=\bold Z\oplus\tau\bold Z$. Up to equivalence, we may suppose that $g_2=g_2(\tau)$ and that $g_3=g_3(\tau)$. Let $\zeta(z)$ and $\eta(z)$ be respectively the Weierstra{\ss} $\zeta$-function and the quasi-period function associated to $\tau$ (cf.~\cite{Silverman}). The function
$$\theta(z)=\zeta(z)-\eta(z)$$
defines a non-holomorphic (but real analytic) map $E\to\P$.

\begin{prop} Setting $K=\bold Q(g_2,g_3)$, the function $\theta$ defines a map
$$E(\CC)_\tors\setminus E[2]\to K$$
Its zeroes correspond to the Lam\'e operators $L_1$ associated to $E$ with dihedral projective monodromy group.
\end{prop}

\begin{proof} Let $\Div^0(E)'$ denote the group of degree $0$ divisors on $E$ which are prime to $0_E$ and consider the map
$$\aligned\Div^0(E)'&\to E\times\Ga\\
\sum_{i}n_i[P_i]&\mapsto\left(\sum_in_iP_i,\sum_in_i\theta(P_i)\right)\endaligned$$
By endowing $E\times\Ga$ with its natural structure of group, the above map is in fact a homomorphism and the general properties of elliptic functions imply that its kernel is precisely the group of principal divisors $D=(t)$ with $t\in 1+\frak m^2$, so that we obtain a real analytic (but not holomorphic, nor algebraic) isomorphism $J_\frak m\cong E\times\Ga$. We know from Proposition~\ref{prop1} that an operator $L_1$ with dihedral projective monodromy group corresponds to a point $P\in E(\CC)_\tors\setminus E[2]$ such that $\varphi(P)$ is a torsion point of $J_\frak m$. By identifying $J_\frak m$ with $E\times\Ga$, we find $\varphi(P)=(2P,2\theta(P))$ and thus $\varphi(P)\in J_{\frak m,\tors}$ if and only if $\theta(P)=0$. The fact that $\theta(P)$ belongs to $K$ for any torsion point $P$ of $E$ is proved in~\cite{Baker} or in~\cite{Mall}.
\end{proof}

\section{Grothendieck's dessins d'enfants}

Considered as a purely combinatorial object, a {\it dessin d'enfant} (litterally, a {\it child's drawing}) is an abstract (connected) graph endowed with two extra structures: a cyclic ordering of the edges meeting at a same vertex and a bipartite structure on the set of its vertices, i.e. a distinction between black and white vertices in such a way that the two ends of any edge never have the same color. Following the ideas exposed by A. Grothendieck in his ``Esquisse d'un programme''~\cite{Grothendieck_1}, these objects classify the isomorphism classes of covers of the projective line (over $\bold C$) which are unramified outside the points $\infty,0$ and $1$. This correspondence is obtained via the topological theory of the fundamental group. The {\it degree} of a dessin d'enfant is the number of its edges, which coincides with the degree of an associated cover.

A rigidity criterion of A. Weil~\cite{Weil} asserts that each isomorphism class of \'etale covers of $\PP=\P_\CC\setminus\{\infty,0,1\}$ has a representative defined over $\overline{\bold Q}$, on which the absolute Galois group $\GQ=\rm{Gal}(\overline{\bold Q}/\bold Q)$ acts in a natural way. Such an action is compatible with the notion of isomorphism and induces a Galois action on the set of dessins d'enfants which translates the action of $\GQ$ on the algebraic fundamental group of $\bold P^*$. It is then possible to introduce the {\it field of definition} (usually called {\it field of moduli}) of a dessin d'enfant, which in most of the cases is the smallest field of definition for the associated covers. Its degree is just the number of Galois conjugates of the dessins d'enfant.

In the following, we call {\it tree} a dessin d'enfant with no closed loops. It corresponds to an isomorphism class of covers $\P_\CC\to\P_\CC$ totally ramified above the point $\infty$ and one can easily prove that its field of moduli is in fact a field of definition. In this paper, we are concerned with the following particular class of trees: given three positive integers $a,b$ and $c$, we denote by $[a,b,c]$ the only tree of degree $N=a+b+c$ having one ``central'' black vertex of valency $3$ and three ``branches'' made of $a,b$ and $c$ edges respectively (turning around the central vertex counterclockwisely, see the following figure). We clearly have $[a,b,c]=[b,c,a]=[c,a,b]$. The {\it signature} of the tree $[a,b,c]$ is the number of its black vertices of valency $1$, its {\it order} is the integer $N/\gcd(a,c,b)$. We say that the tree is {\it primitive} if its degree is equal to its order.

\vskip.4cm
\begin{center}
\includegraphics[scale=.8]{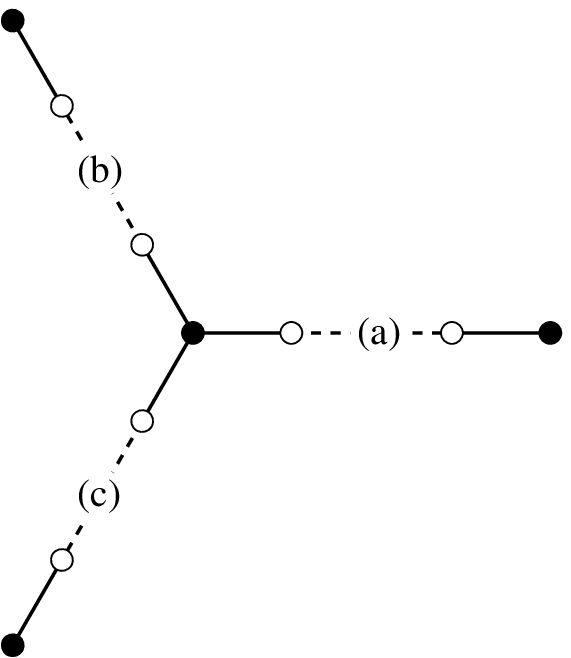}\\
{Figure 1. The tree $[a,b,c]$}
\end{center}
\vskip.4cm

The degree and the signature are clearly Galois invariants. The fact that the order is also invariant under the action of $\GQ$ is less trivial, see for example~\cite{Pakovitch} or~\cite{Zapponi_1}. One can moreover easily check that the complex conjugation sends the tree $T=[a,b,c]$ to the tree $[a,c,b]$, so that $T$ is defined over $\bold R$ if and only if at least two of the integers $a,b$ and $c$ are equal.

\begin{theo}\label{theo1} For any positive integer $N$, there is a one-to-one correspondence between the set of equivalence classes of Lam\'e operators $L_1$ with dihedral projective monodromy group of order $2N$ and the set of primitive trees $[a,b,c]$ with $a+b+c=N$. Moreover, the field of moduli of such an operator coincides with the field of moduli of the corresponding tree.
\end{theo}

\begin{proof} We know from Lemma~\ref{lemm1} and Proposition~\ref{prop1} that an equivalence class of Lam\'e operators $L_1$ with dihedral projective monodromy of order $2N$ corresponds to an isomorphism class of pairs $(E,P)$ satisfying the second condition of Proposition~\ref{prop1}, their field of moduli being the same. Now, the results in~\cite{Zapponi_1} assert that there is a one-to-one correspondence between the set of isomorphism classes of such pairs and the primitive trees $[a,b,c]$ of degree $N$; once again the fields of moduli coincide.
\end{proof}

From a practical and explicit point of view, the correspondence of Theorem~\ref{theo1} can be obtained as follows: consider a primitive tree $[a,b,c]$ of degree $N$ and let $\beta:\P\to\P$ be a model associated to it. Since the cover is totally ramified above $\infty$, we may assume that it is induced by a polynomial $\beta(x)\in\CC[x]$. Let $S\subset\P(\CC)$ be the set of elements of $\beta^{-1}(\{0,1\})$ with odd ramification index (the center and the three ends of the tree, cf. Figure 1). The elliptic curve $E$ is realized as the unique (up to isomorphism) double cover $\pi:E\to\P$ having $S$ as branch locus. Its origin $0_E$ is, by definition, the preimage under $\pi$ of the center of the tree and, more generally, we find $E[2]=\pi^{-1}(S)$. As before, we denote by $\sigma$ the canonical involution of $E$. Since $\infty$ does not belong to $S$, we have $\pi^{-1}(\infty)=\{P,\sigma(P)\}$ and one checks that $P$ (or $\sigma(P)$) satisfies the second condition of Proposition~\ref{prop1}. Conversely, let $L_1=L_{1,g_2,g_3,B}$ be a Lam\'e operator with dihedral projective monodromy of order $2N$ and fix an element $P$ such that $B=x(P)$, so that it satisfies the second condition of Proposition~\ref{prop1}. Let $t$ be the unique function such that $(t)=N[P]-N[\sigma(P)]$ and $v_{0_E}(t-1)=3$. The induced cover $E\to\P$ is unramified outside $\infty,0$ and $1$. Since $\sigma^*t=t^{-1}$, we deduce that the rational function $\beta_0=-(t-1)^2/4t$ is invariant under $\sigma$ and thus, we obtain a commutative diagram
$$\xymatrix{
& E\ar[r]^t\ar[d]_\pi\ar[rd]^{\beta_0}&\P\ar[d]^{-{(x-1)^2\over4x}}\\
& \P\ar[r]^\beta &\P
}$$
It then follows from Abhyankar's Lemma that the cover $\beta$ is unramified outside the set $\{\infty,0,1\}$ and that it is a model for a primitive tree $[a,b,c]$ of degree $N$. For further details, see~\cite{Zapponi_1}.

\begin{coro}\label{cor2} For any positive integer $N$, there are finitely many equivalence classes of Lam\'e operators $L_1$ with dihedral projective monodromy of order $2N$.
\end{coro}

\begin{proof} This follows from the fact that there exists finitely many primitives trees $[a,b,c]$ with $a+b+c=N$.
\end{proof}

\begin{coro}\label{cor1} The field of moduli of a Lam\'e operator $L_1$ with dihedral projective monodromy group is a number field. In particular, there is a natural action of ${\rm Gal}(\QQ/\bold Q)$ on the set of equivalence classes of such operators.
\end{coro}

\begin{proof} Indeed, the field of moduli of any dessin d'enfant is a number field. The Galois action on the equivalence classes follows from the Galois action on trees.
\end{proof}

The correspondence of Theorem~\ref{theo1} can be generalized to all Lam\'e operators with finite monodromy, allowing a direct enumeration of them. This strategy was succesfully adopted by R. Li\c tcanu in~\cite{Litcanu_1,Litcanu_2}, see also the recent work of S. Dahmen~\cite{Dahmen}.

\section{Some general properties of the field of moduli}

As we have seen in Corollary~\ref{cor2}, up to equivalence, there exists finitely many Lam\'e operators $L_1$ with fixed dihedral projective monodromy. We start this section by giving a similar finiteness result in terms fo the degree of the field of moduli.

\begin{prop}\label{prop3} For any positive integer $d$, there exist finitely many equivalence classes of Lam\'e operators $L_1$ with dihedral projective monodromy having a field of moduli of degree less than or equal to $d$. 
\end{prop}

\begin{proof} We know from a theorem of Merel~\cite{Merel} that there exists a constant $C=C(d)$ only depending on the integer $d$ such that for any number field $K$ of degree less than or equal to $d$ and for any elliptic curve $E$ defined over $K$, the cardinality of $E(K)_\tors$ is bounded by $C$. Suppose now that $L_1$ is an operator with dihedral projective monodromy of order $2N$ and that its field of moduli $K$ has degree less than or equal to $d$. Up to equivalence, we can suppose that the curve $E$ and the element $B=x(P)$ are defined over $K$. In particular, the point $P$ is defined over a number field of degree $d'\leq 2d$ and the same holds for the point $Q=2P$, which has exact order $N$ (cf. the proof of Proposition~\ref{prop1}). This implies that the integer $N$ is bounded by a constant only depending on $d$ and the proposition follows from Corollary~\ref{cor2}.
\end{proof}

Recall that an elliptic curve defined over a number field $K$ has {\it potentially good reduction} at a prime $\frak p$ of $\mathcal O_K$ if its $j$-invariant belongs to the localization of $\mathcal O_K$ at $\frak p$. This means that there exists a finite extension $L/K$ and a model of $E$ over $\mathcal O_L$ which has good reduction at any prime $\frak q$ lying above $\frak p$. An elliptic curve defined over a perfect field $k$ of characteristic $p>0$ is {\it ordinary} if its full $p$-torsion subgroup is non-trivial (and thus cyclic of order $p$), otherwise it is {\it supersingular}. The ordinarity of the curve only depends on its $j$-invariant (cf.~\cite{Silverman_1}). We say that an elliptic curve $E$ defined over $K$ has {\it potentially ordinary reduction} (resp. {\it potentially supersingular reduction}) at $\frak p$ if it has potentially good reduction at $\frak p$ and if there exists an integral model (defined over the ring of integers of a finite extension of $K$) of the curve with ordinary (resp. supersingular) reduction at a prime above $\frak p$. These notions only depend on the image of the $j$-invariant of $E$ in the residue field of $\frak p$ and not on the given model. With a slight abuse of language, the curve has {\it good} ({\it ordinary} or {\it supersingular}) reduction at $\frak p$ if there exists a model $\mathscr E/\mathcal O_K$ of $E$ which has good (ordinary or supersingular) reduction at $\frak p$\footnote{By model, we mean a proper flat scheme $\mathscr E$ over $\mathcal O_K$ for which the generic fiber is only $\QQ$-isomorphic to $E$ and not $K$-isomorphic, as it is usually the case.}.

\begin{prop}\label{prop4} Let $\frak p$ be a prime of the field of moduli $K$ of a Lam\'e operator $L_1$ with dihedral projective monodromy of order $2N$ and denote by $p$ its residual characteristic. If $p>N$ then the extension $K/\bold Q$ is unramified at $\frak p$ and the curve $E$ has good reduction at $\frak p$.
\end{prop}

\begin{proof} First of all, after completion, we can reduce to the case where $K$ is a $p$-adic field.  We denote by $R$ its ring of integers and by $k=R/\frak p$ its residue field. Suppose that $L_1$ corresponds to a pair $(E,P)$ and let $t\in K(E)$ be the associated rational function (cf. the proof of Proposition~\ref{prop1}). The induced cover $t:E\to\P_K$ is of degree $N$ and unramified outside $\infty,0$ and $1$. Its monodromy group can be realized as a subgroup of the symetric group $S_N$ and since $p>N$, we deduce that its order if prime to $p$. In this case, the results of~\cite{Beckmann} (see also~\cite{Chambert-Loir,Emsalem,Wewers}) assert that the cover has good reduction at $\frak p$, i.e. that there exists a smooth $R$-model $E_R\to\P_R$ of the cover $t$ and from this we classically deduce that the field of moduli is unramified at $\frak p$ (cf. {\it loc. cit.}). 
\end{proof}

\begin{rem} For general $n\in\bold Z$, the same arguments show that the field of moduli of a Lam\'e operator $L_n$ with dihedral projective monodromy of order $2N$ is unramified outside the primes which are less than or equal to $nN$.
\end{rem}

\section{The case of dihedral projective monodromy of order $2p$}

We now restrict to the case of Lam\'e operators $L_1$ with dihedral projective monodromy of order $2p$, where $p$ is an odd prime number (we already know from~\cite{Baldassarri} that the case $p=2$ is impossible, this also follows from Theorem~\ref{theo1}, since a primitive tree $[a,b,c]$ has degree at least $3$). Remark that the identity $a+b+c=p$ is sufficient for ensuring that the tree $[a,b,c]$ is primitive. Moreover, its signature is equal either to $0$ or to $2$. In particular, by using the correspondence of Theorem~\ref{theo1}, one can easily show that for $p>3$ there are exactly $(p-1)(p-2)/6$ equivalence classes of such operators; $(p^2-1)/24$ of them correspond to trees with signature $0$ while the remaining $(p-1)(p-3)/8$ are associated to trees with signature $2$. The first result of this section gives a lower bound for the ramification index of a prime in the field of moduli lying above $p$.

\begin{theo}\label{theo2} Let $L_1$ be a Lam\'e operator with dihedral projective monodromy of order $2p$, with $p>3$ prime. Fix a prime $\frak p|p$ of its field of moduli $K$, denote by $e_\frak p$ its absolute ramification index and set
$$e=\frac{p+1-s}{\gcd(p+1-s,4(3-s))}$$
where $s$ is the signature of the tree associated to $L_1$. Then the integer $e$ divides $e_\frak p$.
\end{theo}

\begin{proof} The results in~\cite{Zapponi_2} assert that, given a tree of prime degree $p$, the integer $e_\frak p$ is a multiple of an integer only depending on the ramification data, which, in this case, coincides with $e$.
\end{proof}

The following table gives the possible values of $e$ depending on the residue class of $p$ modulo $12$ and on the signature of the tree associated to the Lam\'e operator.
\vskip.4cm
\begin{center}
\begin{tabular}{|c|c|c|}
\hline
\rule[5mm]{0mm}{-1mm}
$p\mod 12$ & Signature & $e$\\[.05cm]
\hline
\rule[-2mm]{0mm}{6mm} $1,9$ & $0$ & $\frac{p+1}2$ \\
\hline
\rule[-2mm]{0mm}{6mm} $1,5,9$ & $2$ & $\frac{p-1}4$\\
\hline
\rule[-2mm]{0mm}{6mm} $3,7$ & $0$ & $\frac{p+1}4$\\
\hline
\rule[-2mm]{0mm}{6mm} $3,7,11$ & $2$ & $\frac{p-1}2$\\
\hline
\rule[-2mm]{0mm}{6mm} $5$ & $0$ & $\frac{p+1}6$\\
\hline
\rule[-2mm]{0mm}{6mm} $11$ & $0$ & $\frac{p+1}{12}$\\
\hline
\end{tabular}
\end{center}
\vskip.25cm

We now investigate the reduction behaviour of the curve $E$.

\begin{theo}\label{theo3} The assumptions and notation being as in Theorem~\ref{theo2}, the curve $E$ always has potentially good reduction at $\frak p$. Moreover, if $P\in E[2p]\setminus E[2]$ denotes the point associated to $L_1$ (cf.~\S1), the following conditions are equivalent:
\begin{enumerate}
\item The curve $E$ has potenitally supersingular reduction at $\frak p$. 
\item The associated tree has signature $0$.
\item The point $P$ has order $p$.
\item The (full) monodromy of $L_1$ coincides with its projective monodromy.
\end{enumerate}
\end{theo}

\begin{proof} As in the proof of Proposition~\ref{prop4}, we may assume that $K$ is a $p$-adic field. Fix  a model $\beta:\P\to\P$ associated to the tree corresponding to $L_1$ and defined over $K$. The theorem is proved by investigation of the stable model of the cover $\beta$. The results in~\cite{Zapponi_3} (which generalize, in an arithmetic-geometric setting, the earlier works in~\cite{Zapponi_2}) asserts that the (special fiber of the) minimal semi-stable model of $\beta$ which separates the elements of the ramified fibers can be described as shown in Figure 2.
\vskip.4cm
\begin{center}
\includegraphics[scale=.6]{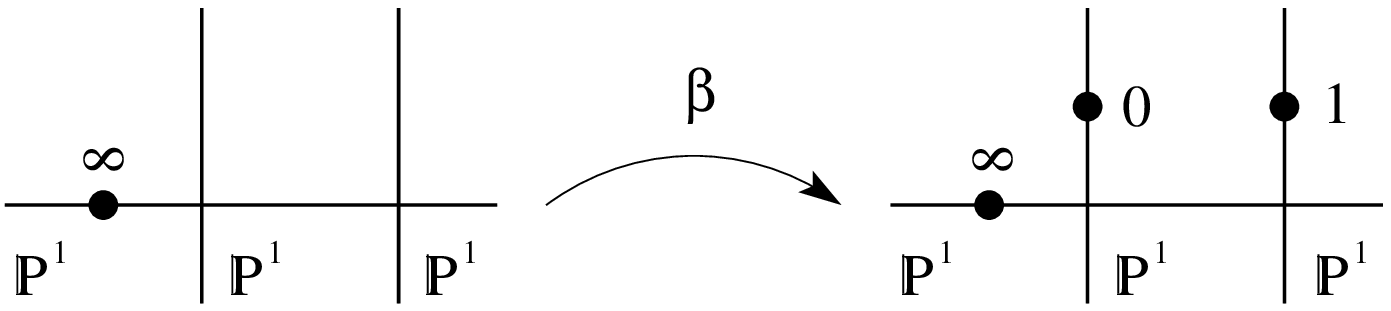}\\
{Figure 2. Semi-stable model for $\beta$}
\end{center}
\vskip.4cm

\noindent More precisely, if $P_0\in\P(\QQ)$ corresponds the center of the tree and if $P_1,P_2$ and $P_3$ are the points associated to its ends then we find the following two possibilities: first of all, if the tree has signature $0$ then $P_0$ lies in the fiber of $\beta$ above $0$ while $P_1,P_2$ and $P_3$ are mapped to $1$. This is the case described in Figure 3.

\vskip.4cm
\begin{center}
\includegraphics[scale=.6]{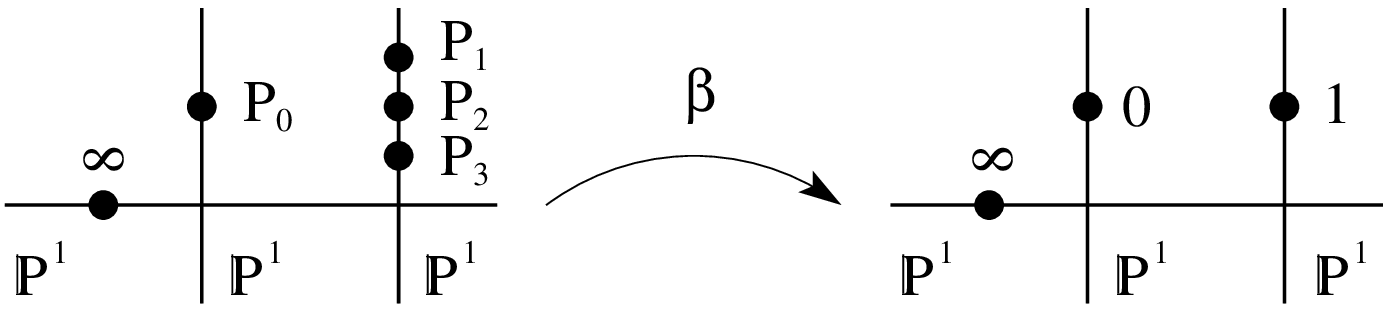}\\
{Figure 3. Semi-stable model in the case of signature $0$}
\end{center}
\vskip.4cm

\noindent Finally, if the tree has signature $2$ then, up to a permutation of the points $P_1,P_2,P_3$, we can assume that $\beta(P_0)=\beta(P_1)=\beta(P_2)=0$ and $\beta(P_3)=1$ and Figure 4 describes the behaviour of the special fiber of the corresponding semi-stable model.

\vskip.4cm
\begin{center}
\includegraphics[scale=.6]{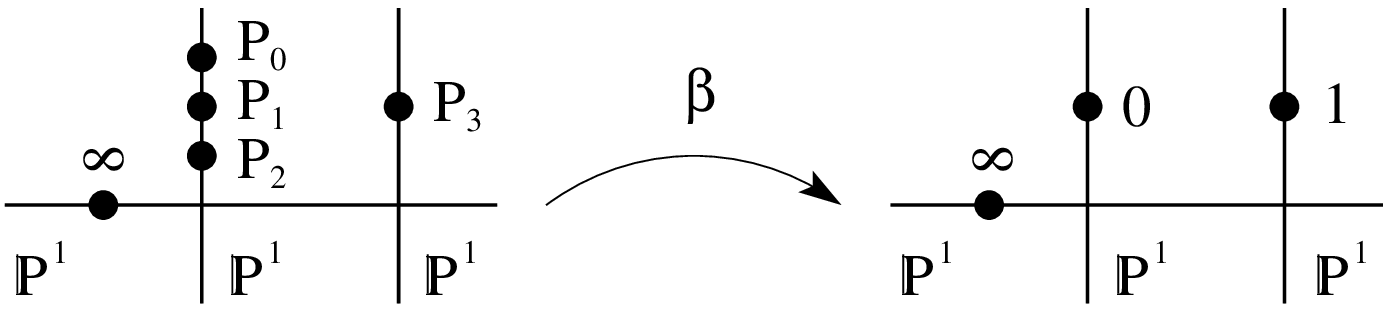}\\
{Figure 4. Semi-stable model in the case of signature $2$}
\end{center}
\vskip.4cm

\noindent In both cases, we see that the points $P_0,P_1,P_2,P_3$ have (potentially) good reduction and since the curve $E$ is realized, up to isomorphism, as the double cover of the projective line ramified at these four points (cf. \S3), we deduce that it has potentially good reduction. This proves the first part of the theorem. The above description for the semi-stable model of $\beta$ allows to deduce the semi-stable model for the cover associated to the function $t$ in the proof of Proposition~\ref{prop1} (see the commutative diagram in \S3). Skipping the details, if the signature of the tree is equal to $0$ then, in such a model, the reduced curve $\overline E$ appears as a degree $p$ cover of the projective line uniquely (and wildly) ramified above one point and the results in~\cite{Zapponi_4} assert that $\overline E$ is supersingular. If the signature is equal to $2$ then the curve $\overline E$ is realized as a degree $p$ cover of the projective line unramified outside two points, wildly ramified above one of them and tamely ramified above the other with a unique (effectively) ramified point over it, with ramification index $3$ and we can easily deduce from the results in {\it loc. cit.} that $\overline E$ is ordinary. This shows the equivalence between the conditions 1 and 2 of the theorem. The equivalence of the conditions 2 and 3 is a restatement of the porposition in \S2.3 of~\cite{Zapponi_1}. Finally, an explicit expression of the solutions of the Lam\'e equation given in~\cite{Beukers} shows that the order of the full monordomy group is the order of the point $P$. 
\end{proof}

\begin{rem} The potentially good reduction of the elliptic curve $E$ can be easily deduced from S. Wewers' results in~\cite{Wewers} which allow to directly determine the semi-stable model of the cover $t$ associated to the torsion point $P$ (cf. \S\S 2,3) without reducing to the genus zero case. The supersingularity criterion follows from a more detailed study of the action of the Cartier operator on the differential forms on $E$ having only one pole.
\end{rem}

\begin{coro}\label{cor3} The potential supersingularity of the reduction of the curve $E$ at $\frak p$ is independent of the prime $\frak p$ of $K$ lying above $p$. 
\end{coro}

\begin{proof} Immediate, since the prime $\frak p$ does not appear in the conditions 2 and 3 in Theorem~\ref{theo3}.
\end{proof}

We know that, up to a finite extension of the field of moduli $L/K$, the curve $E$ has good reduction at any prime $\frak p|p$ of $L$. On the other hand, it also admits a model over $K$, for which the good reduction at $\frak p$ is not ensured. The following result gives some furhter informations and relates the behaviour of the reduction to the ramification index $e_\frak p$ at $\frak p$.

\begin{prop}\label{prop5} Let $L_1$ be a Lam\'e operator with dihedral projective monodromy of order $2p$ associated to an elliptic curve $E$ with invariant $j$. Set $m=2ns$, where $s$ is the signature of the correspnding tree and $n=3,2$ if $j=0,1728$ and $n=1$ otherwise. Let $\frak p$ be a prime of the field of moduli $K$ lying above $p$. The following conditions are equivalent:
\begin{itemize}
\item The curve $E$ has good reduction at $\frak p$.
\item The integer $(p+1-s)/m$ divides $e_\frak p$.
\end{itemize}
\end{prop}

\begin{proof} As usual, we may suppose that $K$ is a $p$-adic field with ring of integers $R$. We moreover fix an algebraic closure $\overline K$ of $K$. Let $|\,\,\,|_\frak p$ be the $\frak p$-adic norm, normalized by the condition $|p|_\frak p=p^{-1}$ and (uniquely) extended to the whole $\overline K$. Let also $v_\frak p$ be the associated valuation, with $v_\frak p(p)=e_\frak p$. Suppose that the signature of the corresponding tree is $0$ and that $j\neq0,1728$. We know from~\cite{Zapponi_1} that there exists a polynomial model $\beta:\P_K\to\P_K$ of it, which can be written as
$$\beta(x)=x^3g(x)^2=1+f(x)h(x)^2$$
with $f,g,h\in R[x]$, $f$ monic of degree $3$ and $g,h$ of degree $(p-3)/2$ (we moreover assume that the leading coefficients of $f$ and $g$ are $R$-units). The point $x=0$ corresponds to the center of the tree while the roots of $f$ are its ends. The main result in~{\it loc. cit.} asserts that, for any two distinct roots $x_1,x_2\in\overline K$ of $\beta-1$, we have
$$|x_1|_\frak p=|x_2|_\frak p=0\quad{\rm and}\quad|x_1-x_2|_\frak p=p^{-\frac2{p+1}}$$
In particular, this holds for the roots $x_1,x_2,x_3$ of $f$. The curve $E$ is obtained as double cover of the projective line ramified at the points $0,x_1,x_2,x_3$. We can replace these four points by their image under the action of an element $\tau$ of ${\rm PGL}_2(K)$, the resulting ellitpic curve will be isomorphic to $E$; taking $\tau=1/z$ we obtain the elements $\infty,y_1,y_2,y_3$, where we have set $y_i=\tau(x_i)$. The curve $E$ is then given by the affine Weierstra{\ss} equation
$$y^2=f_1(x)=(x-y_1)(x-y_2)(x-y_3)$$
By construction, we have $|y_i|_\frak p=0$ and $|y_i-y_j|_\frak p=p^{-{2\over p+1}}$ for $i\neq j$. This implies that the discriminant $\Delta$ of $f_1$ satisfies the identity
$$|\Delta|_\frak p=p^{-{12\over p+1}}$$
which can be rewritten as $v_\frak p(\Delta)=12e_\frak p/(p+1)$. Now, we know from~\cite{Silverman_1} (see also~\cite{Silverman}) that there exists a smooth $R$-model of $E$ if and only if $v_\frak p(\Delta)\equiv0\pmod6$\footnote{Recall that we are concerned with $\overline K$-isomorphism classes of elliptic curves and not with $K$-isomorphism classes. That's why we consider the valuation of the discriminant modulo $6$ instead of $12$, allowing quadratic twists of $E$.}. This last identity is equivalent to $e_\frak p\equiv 0\pmod{\frac{p+1}2}$, as desired. The case of signature $2$, as the cases $j=0$ and $j=1728$ are treated similarly.
\end{proof}

\begin{exa} We close the paper with the description of the equivalence classes of Lam\'e operators with dihedral projective monodromy group of order $14$. There are $5$ of them, corresponding to the trees $[1,1,5], [1,3,3]$ (signature $0$, both of them are defined over $\bold R$), $[1,2,4],[1,4,2]$ and $[2,2,3]$ (signature $2$, only the last is defined over $\bold R$). In the case of signature $0$, Theorem~\ref{theo2} asserts that the ramification index of any prime above $7$ in the field of moduli is divisible by $8/\gcd(8,12)=2$. Since there are exactly two such trees and since they are not defined over $\bold Q$, we deduce that they are Galois conjugates and that the field of moduli $K$ is totally ramified above $7$ (a direct calculation gives $K=\bold Q(\sqrt{21})$). In this case, we obtain $(p+1-s)/2=4$ and thus Theorem~\ref{theo3} implies that the associated elliptic curves don't have good reduction at the prime of $K$ lying above $7$. Nevertheless, they admit a Weierstra{\ss} model over $K$ with discriminant of valuation $3$ (cf. the proof of Theorem~\ref{theo3}) and Tate's algorithm~\cite{Silverman} asserts that the corresponding N\'eron models are of type III. Remark that there may also exist models with disciminant $9$ having N\'eron models of type III$^*$, but they are obtained as twists of the previous ones. Let's now study the case of signature $2$: the ramification index of any prime above $7$ is divisible by $6/\gcd(6,4)=3$ and thus, since there are three trees of this type, we deduce that they form a unique Galois orbit and that their field of moduli are totally ramified above $7$. Finally, we have $(p+1-s)/2=3$ and thus the curves have good reduction at the unique prime above $7$. Remark that, once again, some twists may have bad reduction, with N\'eron model of type I$_0^*$.
\end{exa}

\end{document}